\documentclass[10pt]{article}
\usepackage{amsfonts,amssymb,latexsym, amsmath, amsthm}

\usepackage{amsfonts,amssymb,latexsym, amsmath, amsthm}

\usepackage[usenames]{color}%%  NEW

\newtheorem{thm}{Theorem}[section]
\newtheorem{lem}[thm]{Lemma}
\newtheorem{cor}[thm]{Corollary}

\theoremstyle{remark}

\def\ZZ{{\mathbb Z}}
\def\FF{{\mathbb F}}

\font\teneufm=eufm10
\font\seveneufm=eufm7
\font\fiveeufm=eufm5
\newfam\eufmfam
\textfont\eufmfam=\teneufm
\scriptfont\eufmfam=\seveneufm
\scriptscriptfont\eufmfam=\fiveeufm

\def\FF{{\mathbb F}}
\def\leq{\leqslant}
\def\geq{\geqslant}
\def\phi{\varphi}
\def\epsilon{\varepsilon}

\begin{document}

\def\normO{{\triangleleft_O\,}}

\def\C{{\rm C}}
\def\Z{{\rm Z}}
\def\cT{{\mathcal T}}
\def\cU{{\mathcal U}}

\def\wr{\,{\rm wr}\,}
\def\L{{\mathcal L}}

\let\oldcomma=,
\catcode`\,=\active
\def,{\ifmmode \oldcomma \else {\/\rm\oldcomma}\fi}
\let\oldcolon=:
\catcode`\:=\active
\def:{\ifmmode \oldcolon \else {\/\rm\oldcolon}\fi}
\let\oldscolon=;
\catcode`\;=\active
\def;{\ifmmode \oldscolon \else {\/\rm\oldscolon}\fi}

\def\Cr{{\rm Cr}}
\centerline{PROFINITE GROUPS WITH FEW CONJUGACY CLASSES} 

\centerline{OF ELEMENTS OF INFINITE ORDER}
\bigskip

\centerline{JOHN S.\ WILSON}

%\date{11 October 2021}
%\maketitle
\bigskip

{\small \noindent ABSTRACT. \
%\begin{abstract}
It is proved that every finitely generated profinite group with fewer than $2^{\aleph_0}$ conjugacy classes of elements of infinite order is finite.  }  %\end{abstract}
\bigskip

\centerline{1.\ INTRODUCTION}
\bigskip

In \cite{jzn}, Jaikin-Zapirain and Nikolov proved that a profinite group with countably many conjugacy classes must be finite.  
A well-known theorem of Zelmanov \cite{efim} asserts that profinite torsion groups are locally finite; that is, their finite subsets generate finite subgroups.  
 Inspired by these results, we use an extension of Zelmanov's theorem (also due to Zelmanov; see Theorem C below) and some ideas from Wilson \cite{fewccls} to prove the following result:

\begin{description}\item[Theorem A] \begin{em}  Let $G$ be a finitely generated profinite group.  If $G$ has fewer than $2^{\aleph_0}$ conjugacy classes of elements of infinite order, then $G$ is finite. \end{em}\end{description}

We shall also prove a local version of this result, concerning $p$-elements for a fixed prime $p$.
We recall that a (generalized) $p$-element of a profinite group $G$ is an element that (topologically) generates a finite $p$-group or a copy of the group $\ZZ_p$ of $p$-adic integers. By a $p^\infty$-{\em element} we understand a $p$-element of infinite order.    The $p$-Sylow subgroups of $G$ are the maximal subgroups consisting of $p$-elements; such subgroups exist and are conjugate (see for example  \cite[Chapter 2]{book}). 

\begin{description}\item[Theorem B] \begin{em}  Let $G$ be profinite group, $p$ an odd prime and $P$ a $p$-Sylow subgroup.  If $G$ has fewer than $2^{\aleph_0}$ conjugacy classes of $p^\infty$-elements, and $P$ is finitely generated, then  $P$ is finite. \end{em}\end{description}

We conjecture that natural extensions of the above theorems hold without the hypotheses of finite generation; we explain briefly why this hypothesis arises.  When trying to prove that a profinite torsion group is locally finite, one can immediately pass to finitely generated subgroups, and results of Zelmanov on finitely generated Lie algebras can be used.  In our proofs, a critical case is concerned with groups having the property that the set of torsion elements is dense, and this property does not pass to finitely generated or indeed procyclic subgroups (as the existence of infinite abstract finitely generated residually finite torsion groups shows).  The difficulty that arises is reminiscent of the open problem of determining whether profinite torsion groups have finite exponent: this too is essentially a problem about (countably based) pro-$p$ groups that are not finitely generated.  

\bigskip

{\hrule
\smallskip


\newpage

\centerline{2. \ PRELIMINARIES; COSETS OF CONJUGATES}\medskip

For unexplained notation and general information about profinite groups and their Sylow theory, we refer the reader to \cite{book}.  We write $N\normO G$ to indicate that $N$ is an open normal subgroup of a profinite group $G$, and for
$x\in G$ we write $x^G$ for the conjugacy class of $x$ and $Nx^G$ for the product of the sets $N$, $x^G$.  Each conjugacy class $x^G$ is closed (being the image of $G$ under the continuous map $g\mapsto x^g$). 
All subgroups arising are understood to be closed.  A subset of a group is said to have {\em finite exponent} if its elements have bounded orders.   The following elementary facts will be used frequently: elements $x,y$ of a profinite group $G$ are conjugate in $G$  if and only if for each $M \normO G$ the elements $Mx,My$ are conjugate in $G/M$ (see e.g.\ \cite[Lemma 2.1]{fewccls});
the property of having fewer than $2^{\aleph_0}$ conjugacy classes of elements of infinite order is inherited by open subgroups and continuous images.

\begin{description} \begin{em}\item[Proposition 2.1.]   Let $G$ be a profinite group, $p$ a prime and $P$ a $p$-Sylow subgroup of $G$.    Let $N_0\normO G$, $u\in P$.

Suppose that, for each coset $Nt$ with $N\normO G$, $t\in P$ and $Nt\subseteq N_0u$

{\rm(a)} $Nt\cap P\not\subseteq t^G$, and

{\rm(b)}   $Nt\cap P$ does not have finite exponent.

\noindent Then $G$ has at least $2^{\aleph_0}$ $G$-conjugacy classes of $p^\infty$-elements.

\end{em}\end{description}
\begin{proof}
The result is very similar to Proposition 2.2 in \cite{fewccls}, which was the special case when $N_0u=G$, and the proof is essentially the same.  
However the printed proof in \cite{fewccls} contains small errors (corrected in the arXiv version) and so we give the proof in its entirety.   

We construct a descending chain $(N_k)_{k\geq0}$ of open normal subgroups and a family $(R_k)_{k\geq0}$ of finite subsets of
$N_0u\cap P$ such that for each $k\geq1$ and $x\in R_k$ there are elements $x^{(1)},x^{(2)}\in N_{k-1}x\cap P$ for which 
\smallskip

\ (i)  $N_k x^{(1)}$ and $N_k x^{(2)}$ are not conjugate in $G/N_k$ and

(ii) $N_kx^{(1)}$ and $N_kx^{(2)}$ have order at least $p^k$ in $G/N_k$.

\medskip

We set $R_0=\{u\}$.  
Suppose that $k\geq1$ and that $N_{k-1}$, $R_{k-1}$ have been constructed, and let $x\in R_{k-1}$.  By (b) we can find $x^{(1)}\in N_{k-1}x\cap P$ and a subgroup $L_x\normO G$ with $L_x< N_{k-1}$ such that $L_xx^{(1)}$ has order at least $p^k$ in 
$G/L_x$.  By (a), the result \cite[Lemma 2.1]{fewccls} mentioned above and again (b), we can find $M_x\normO G$ with $M_x\leq L_x$  and $x^{(2)} \in L_xx\cap P$ such that $M_xx^{(1)}$ and $M_xx^{(2)}$ are non-conjugate elements of $G/M_x$ and with $M_xx^{(2)}$ of order at least $p^k$ in 
$G/M_x$.  We define $$N_k=\textstyle{\bigcap_{x\in R_{k-1}}}M_x\quad\hbox{and}\quad R_k=\{x^{(1)},x^{(2)}\mid x\in R_{k-1}\}.$$  

Now consider the set $F=\{ N_kx^G\mid k\geq0, x\in R_k\}$, partially ordered with respect to inclusion.  By construction, each element of $F$ contains at least two maximal elements; it follows that $F$ has
$2^{\aleph_0}$ maximal chains $(N_kx_k^G)_{k\geq0}$.  By compactness, the intersection $\bigcap N_kx_k^G$ of each such chain is non-empty; it is also evidently a union of conjugacy classes.  
If $(N_kx_k^G)$,  $(N_ky_k^G)$ are distinct chains then their intersections  $\bigcap N_kx_k^G$,  $\bigcap N_ky_k^G$ are disjoint: if $k$ is minimal with
 $N_kx_k^G\neq N_ky_k^G$, then $k\geq1$ and $x_k,y_k$ are distinct elements of $R_k$; since $N_{k-1}x_k^G=N_{k-1}y_k^G$ their relationship is that of the elements $x^{(1)},x^{(2)}$  appearing in the construction of $R_k$ and so $N_kx_k^G$ and $N_ky_k^G$ are disjoint. 
 
 Finally, if $x\in\bigcap_k N_kx_k^G$ then for each $r\geq 1$ we have $N_{r+1}x^{p^r}= N_{r+1}x_{r+1}^{p^r}\neq N_{r+1}$ and
 so $x^{p^r}\notin N_{r+1}$.  Therefore the elements of the intersections of chains $(N_kx_k^G)$ are $p^\infty$-elements and the result follows.
\end{proof}

\begin{description} \item[Lemma 2.2.]\begin{em}  Let $G$ be a  profinite group, $p$ a prime and $P$ a $p$-Sylow subgroup. Let $K$ be a closed normal subgroup and $t\in P$, %, $N\triangleleft_OG$ 
and suppose that
$Kt\cap P\subseteq t^G$.
\begin{enumerate}
\item[\rm(a)]  If $p\neq2$ then $K$ is an extension of a pro-$p'$-group by a pro-$p$-group.
\item[\rm(b)]  If $K$ is open and is an extension of a pro-$p'$-group by a pro-$p$-group then 
$t$ has finite order.  
\end{enumerate}
 \end{em}\end{description}

 \begin{proof}  (a) First suppose that $G$ is finite.  We use results proved in \cite{fewccls}.  From \cite[Lemmas 3.3, 3.7]{fewccls} applied to quotient groups $G/K_1$ with $K_1\triangleleft G$ and $K_1\leq K$, we conclude that $K$ has no perfect composition factors of order divisible by $p$, and then the conclusion follows from  \cite[Lemmas 3.3, 4.1]{fewccls}.  (It is worth mentioning that \cite[Lemma 3.7]{fewccls} depends on the classification of the finite simple groups; in particular this is the reason why the prime $2$ is excluded.) 
 
 Now we consider the general case.   For each $N\triangleleft_OG$ let $R_N/N$  be the smallest normal subgroup with $G/R_N$ a $p$-group;  from above $R_N/N$ is a $p'$-group, and so $\bigcap_NR_N$ is a pro-$p'$ group while $G/\bigcap R_N$ is a pro-$p$ group.  
 
 (b)  To prove that $t$ has finite order we may replace $G$ by $G/R$ where $R$ is the smallest normal subgroup of $K$ with $K/R$ a pro-$p$ group, since this quotient group inherits the hypotheses on $G$.  We have now
 $K\leq P$.  %Write $\kappa=|G\colon P|$.   
 Let $L\normO G$ with $L\leq K$ and write 
$\bar G$, $\bar K$, $\bar t$ for the images of $G, K$ and $ t$ in $G/L$.  Since $\bar K\bar t\subseteq \bar t^{\bar G}$ we have
$$|\bar K|=|\bar K\bar t|\leq |\bar G|/|\C_{\bar G}(\bar t)|,$$
and so $$|\C_{G/L}(Lt)|\leq |\bar G\colon \bar K|= |G\colon K|.$$
It follows that $|DL/L|\leq |G\colon K|$ where $D=\C_G(t)$.  This holds for all $L$, and since $\bigcap L=1$ we conclude that $D$ is finite.  But $\langle t\rangle \leq D$.
\end{proof}

 \begin{description}\item[Corollary 2.3.] \begin{em} Let $G$ be profinite, $p$ a prime and $P$ a $p$-Sylow subgroup.  Suppose that $G$ has fewer than $2^{\aleph_o}$ conjugacy classes of $p^\infty$-elements.  If either $p$ is odd or $G$ is a pro-$2$ group, then for all $N_0\normO G$ and $u\in P$ there exist $N\normO G$, $t\in P$ with $Nt\subseteq N_0u$ and 
 $Nt\cap P$ of finite exponent.  \end{em}\end{description}
 
\begin{proof}  By Proposition 2.1 there exist  $N\normO G$, $t\in P$ with $Nt\subseteq N_0u$ and either
$Nt\cap P\subseteq t^G$ or $Nt\cap P$ of finite exponent.  In the former case, $t$ has finite order by Lemma 2.2, so $Nt\cap P$ again has finite exponent.
\end{proof}  

\bigskip\newpage

\centerline{3.\ DENSE TORSION: ZELMANOV'S THEOREMS}

\medskip

Because the open cosets of a profinite group $G$  constitute a base of neighbourhoods,
the set of torsion elements of $G$ is dense if and only if for each $N_0\normO G$ and $u\in P$ the coset $N_0u$ contains an element of finite order.  In particular, the groups studied in Corollary 2.3 have this property.  

In \cite{efim2017} Zelmanov proved a theorem that gives sufficient conditions for certain finitely generated Lie algebras $L$ to be nilpotent.  We fix a finite generating set $X$ and consider the set $X^*$ consisting of $X$ and all iterated products of elements of $X$.  An element $a\in L$ is called {\rm ad}-nilpotent if the linear map ${\rm ad}_La\colon y\mapsto ya$ from $L$ to $L$ is nilpotent.  Zelmanov's theorem (\cite[Theorem 1.1]{efim2017}) asserts that if 
(i) $L$ satisfies a polynomial identity and (ii)
every element of $X^*$ is {\rm ad}-nilpotent, then $L$ is nilpotent. 

The following theorem is a consequence of this deep result.   It is also due to Zelmanov, and we thank him for allowing us to include the proof.

\begin{description} \item[Theorem C.] \begin{em} 
Let $G$ be a finitely generated pro-$p$ group and suppose that

{\rm(i)} the set of torsion elements is dense in $G$, and

{\rm(ii)} $G$ has an open normal subgroup $N$ and an element $a$ such that $Na$ has finite exponent.

Then $G$ is finite.
\end{em}\end{description}

\begin{proof} Consider the Zassenhaus filtration $(G_k)$ of $G$; this is the descending chain of normal subgroups defined by $G_k=G\cap(1+J^k)$, where $J$ is the augmentation ideal of the group algebra $\FF_pG$.  
Write $\L_p(G)=\bigoplus_{k\geq1} L_k$ where $L_k=G_k/G_{k+1}$ for each $k$.  Then $\L_p(G)$ is a Lie algebra over $\FF_p$, with multiplication $L_i\times L_j\to L_{i+j}$ on the direct summands induced by the commutator map in $G$ %and extended by linearity 
(see for example \cite {efim2017} or \cite{EfimandI}).  

By the main theorem in \cite{EfimandI} and condition (ii), the Lie algebra $\L_p(G)$ satisfies a polynomial identity.  Let $G_{k+1}u$ be an arbitrary element of $L_k$.  By (i), there is a torsion element $x\in G_{k+1}u$.   Thus $G_{k+1}u=G_{k+1}x$, and  an easy calculation shows that if $x^{p^n}=1$ then $({\rm ad}(G_{k+1}u))^{p^n}=({\rm ad}(G_{k+1}x))^{p^n}=0$.  It follows that all elements in the subspaces $L_k$ are ad-nilpotent.  From the result of \cite{efim2017} stated above we conclude that $\L_p(G)$ is nilpotent, and so $G_k=G_{k+1}$ for some $k$.  It follows from a theorem of Lazard see (\cite[Section 11]{DDMS}) that $G$ is $p$-adic analytic, and so $G$ has an open torsion-free normal subgroup $N$.  However if $L\triangleleft_OG$ and $L<N$ then no coset of
$L$ in $N$ apart from $L$ itself has torsion elements.  Hence $N=1$ and $G$ is finite. 
\end{proof}

The case of finitely generated pro-$p$ groups with few conjugacy classes of $p^\infty$-elements is now easily settled. 

\begin{description} \item[Corollary 3.1.]  \begin{em} Suppose that $G$ is a pro-$p$ group with fewer than $2^{\aleph_0}$ conjugacy classes of $p^\infty$-elements.  Then $G$ is finite.
 \end{em}\end{description}
 
 \begin{proof}  
 From Corollary 2.3 (with $G=P$) it is enough to consider the case when (a) $G$ has a coset $N_0u$ of finite exponent and (b) the set of torsion elements is dense in $G$.  Therefore the result follows immediately from Theorem C.
\end{proof}

\bigskip

\centerline{4. PROOF OF THEOREMS A AND B}

\medskip

We need a slight variant of results in \cite{oldcompact}, noted in \cite[Lemma 2.5]{fewccls}.

\begin{description}\item[Lemma 4.1.] \begin{em}   Let $G$ be a profinite group, $p$ be a prime and $P$ a $p$-Sylow subgroup of $G$.   Let $N\normO G$, $t\in P$  and suppose that $Nt\cap P$ has finite exponent.  Then  
$G$ has a finite series of closed characteristic subgroups in which each factor is one of the following: 
{\rm(i)} a pro-$p$-group; {\rm(ii)} a pro-$p'$-group, or {\rm(iii)} a Cartesian product of isomorphic finite simple groups of order divisible by $p$.  
\end{em}\end{description}

It is well known that Cartesian powers of finite groups are locally finite.

\medskip
\noindent{\em Proof of Theorem $B$.}
%\begin[Proof of Theorem B]{proof}  \end{proof}
%\end{document}
From Corollary 2.3 and Lemma 4.1, the group $G$ has a series of the kind described in Lemma 4.1.  
Since the hypotheses are inherited by quotient groups, induction on the length of such a series allows us to assume that $G$ has a normal subgroup $K$ of one of the types (i), (ii) or (iii) in Lemma 4.1 with the property that the Sylow subgroup $PK/K$ of $G/K$ is finite.  Thus $P\cap K$ is finitely generated, and it will suffice to prove that it is finite.  This clearly holds if $K$ is a pro-$p'$-group, or if $K$ is of type (iii) above since then $K$ is locally finite.   

Suppose then that $K$ is a finitely generated pro-$p$ group with $P/K$ finite. It follows that
$G$ has an open normal subgroup $M$ with $K\leq M\leq N$ such that $M/K$ and $P/K$ are disjoint; the subgroup $M$ inherits the hypothesis on conjugacy classes and $M/K$ is a pro-$p'$-group. 
By the Schur--Zassenhaus theorem for profinite groups (cf.\ \cite[Proposition 2.33]{book}),  we have $M=K\rtimes Q$ for a pro-$p'$ subgroup $Q$.  Write $\Phi$ for the Frattini subgroup of $K$.  Since $K$ is finitely generated, $K/\Phi$ is finite, and an open subgroup $Q_0$ of $Q$ acts trivially on it.  It follows (for example from \cite[Theorem 5.1.4]{gor}) that $Q_0$ acts trivially on $K/L$ for each $L\triangleleft_O K$ with $L\triangleleft M$, and hence acts trivially on $K$.  
The subgroup $C=\C_Q(K)$ is normal in $M$  and since $Q_0\leq C$ the image of $Q$ in $M/C$ is finite.  Hence $KC/C$ is an open pro-$p$ subgroup of $M/C$, and so it inherits the hypothesis on conjugacy classes; therefore it is finite by Corollary 3.1.  Therefore the $p$-Sylow subgroups of both $M$ and $G$ are finite, as required. 
\hfill\qed

\medskip

Now we need a consequence of an important result of Herfort \cite{Wolfgang}.
  
\begin{description}\item[Lemma 4.2.]  \begin{em}  If $G$ is a profinite group with fewer than $2^{\aleph_0}$ conjugacy classes of elements of infinite order, then $G$ has non-trivial $p$-elements for only finitely many primes $p$. \end{em}\end{description}

\begin{proof}  
%Let $G$ be a profinite group having non-trivial $p$-elements for infinitely many primes $p$.  
Herfort's theorem states that a profinite group with non-trivial $p$-elements for infinitely many primes $p$ 
has a procyclic subgroup $A$ with the same property; so $A$ is isomorphic to a Cartesian product $A^*=\Cr_{p\in I}A_p$ of non-trivial procyclic $p$-groups $A_p$ for an infinite set $I$ of primes (see \cite[Chapter 2]{book}).  Elements of $A$ whose images in $A^*$ have distinct infinite supports generate non-isomorphic infinite procyclic subgroups and cannot be conjugate.   However an infinite set has at least $2^{\aleph_0}$ infinite subsets. 
\end{proof}

\medskip
\noindent{\em Proof of Theorem $A$.}
By Lemma 4.2, the set $\pi(G)$ of primes $p$ for which $G$ has non-trivial $p$-elements is finite, and by Corollary 3.1 we can assume that $|\pi(G)|>1$.  Assume the result known for all groups $H$ satisfying the hypothesis and with $|\pi(H)|<|\pi(G)|$.  Choose an odd prime $p\in\pi(G)$.  By Corollary 2.3, $G$ has a finite series as described in Lemma 4.1, and by induction on the length of such a series we can suppose that $G$ has an open normal subgroup $K$ with one of the three types described in Lemma 4.1; this subgroup $K$ is finitely generated and it too has fewer than $2^{\aleph_0}$ conjugacy classes of elements of infinite order.  Hence $K$ is finite: if $K$ is a pro-$p$ group or a pro-$p'$ group this follows since $|\pi(K)|<|\pi(G)|$ and in the remaining case $K$ is finitely generated and locally finite.  \hfill \qed

\medskip

Mathematisches Institut, Universit\"at Leipzig, 04109 Leipzig, Deutschland

\smallskip

\quad and

\smallskip

Christ's College, Cambridge CB2 3BU, United Kingdom

\medskip

E-Mail addresses:  {\tt wilson@math.uni-leipzig.de} and {\tt jsw13@cam.ac.uk}


\begin{thebibliography}{Wie}

%\bibitem{atlas} J.H. Conway, R.T. Curtis, S.P. Norton, R.A. Parker and R.A. Wilson.  {\em Atlas of Finite Groups} (Clarendon Press, Oxford, 1985).

\bibitem{DDMS}  J. D. Dixon, M. P. F. du Sautoy, A. Mann, and D. Segal. {\em Analytic pro-$p$ groups,} \ 2nd edn. \ Cambridge Studies in Advanced Mathematics, vol.\ 61 (Cambridge University Press, 1999).

%\bibitem{frob} G.\ Frobenius.  \"Uber einen Fundamentalsatz der Gruppentheorie.  {\em Berliner  Sitzungsber.} (1903),  987--991.

\bibitem{gor} D.\ Gorenstein. {\em Finite Groups,} 2nd edn. (AMS Chelsea Publ., 1980).

%\bibitem{lyons} D.\ Gorenstein, R.\ Lyons and R.\ Solomon. {\em The Classification of the Finite Simple Groups.}  Amer.\ Math.\ Soc.\ Surveys and Monographs 40 (3) (1998).

%\bibitem{Hall}  P.\ Hall. On a theorem of Frobenius.  {\em Proc.\ London Math.\ Soc.} (2) {\bf 40} (1935),  468--501.
 
\bibitem{Wolfgang} W.\ Herfort. An arithmetic property of profinite groups.  {\em Manuscripta Math.} {\bf37} (1982), 11--17.
  

\bibitem{jzn} A.\ Jaikin-Zapirain and N.\ Nikolov.  An infinite compact Hausdorff group has uncountably many conjugacy classes.  {\em Proc.\ Amer.\ Math.\ Soc.} {\bf147} (2019), 4083--4089.

%\bibitem{lumba} Alexander Lubotzky and Avinoam Mann.  Powerful $p$-groups.  I. Finite groups.  {\em J.\ Algebra} {\bf105} (1987), 484--505.

%\bibitem{shalev} A. Shalev.  On almost fixed point free automorphisms.  {\em J. Algebra} {\bf157} (1993), 271--282.  

\bibitem{oldcompact} John S.\ Wilson.
On the structure of compact torsion groups.  {\em Monatsh.\ Math. \bf96}
(1983), 57--66.

\bibitem{book}   John S.\ Wilson.  {\em Profinite
Groups.} \  London Math.\ Soc.\ Monographs, New Series 19
(Clarendon Press, Oxford, 1998). 

%\bibitem{prob} John S.\ Wilson.  The probability of generating a soluble subgroup of a finite group.  {\em J. London Math. Soc.} (2) {\bf75} (2007), 431--446.

\bibitem{fewccls} John S.\ Wilson.   Profinite groups with few conjugacy classes of $p$-elements.   {\em Proc.\ Amer.\ Math.\ Soc.} {\bf 150} (2022), 3297--3305; arXiv:2204.09936.
%(to appear).

\bibitem{EfimandI} John S. Wilson and E.I.\ Zelmanov. Identities for Lie algebras of pro-$p$ groups
{\em J. Pure Appl. Algebra} {\bf 81} (l992), 103--109.


\bibitem{efim} E.I.\ Zelmanov.  On periodic compact groups. {\em Israel J.\ Math.}  {\bf77} (1992), 83--95.

\bibitem{efim2017}  Efim \ Zelmanov.  Lie algebras and torsion groups with identity. 
{\em J. Comb. Algebra} {\bf1} (2017), 289--340. 
 
 
\end{thebibliography}
\end{document}